\documentclass[12pt]{article}
\usepackage{amssymb,amsmath,amsfonts,t1enc}
\newcommand{\Q}{\mathbb Q}
\newcommand{\R}{\mathbb R}
\newcommand{\C}{\mathbb C}
\newcommand{\Z}{\mathbb Z}
\newcommand{\F}{\textbf{\textit{F}}}
\newcommand{\K}{\textbf{\textit{K}}}
\newcommand{\LL}{\textbf{\textit{L}}}
\newcommand{\wt}{\widetilde}

\begin{document}
\thispagestyle{empty}
\hfill
\vskip 1.5truecm
\footnotetext{
{\bf Key words and phrases:} algebraic function field in one
variable over $\Q$ (over $\R$, over ${\Q}_p$),
element transcendental over $\Q$ (over~$\R$, over ${\Q}_p$),
existentially \mbox{$\emptyset$-definable} element,
Faltings' finiteness theorem, field finitely generated \mbox{over $\Q$,}
Pythagorean subfield of $\R$,
recursively approximable real number,
subset \mbox{of $\R$} which is implicitly \mbox{$\emptyset$-definable}
\mbox{in $(\R,+,\cdot,0,1)$.}}
\par
\noindent
\centerline{{\large On $\emptyset$-definable elements in a field}}
\vskip 1.0truecm
\centerline{{\large Apoloniusz Tyszka}}
\vskip 1.0truecm
\par
\noindent
{\bf Abstract.} We develop an arithmetic characterization of
elements in a field which are first-order definable by a parameter-free
existential formula in the language of rings. As applications we show
that in fields containing an algebraically closed field only the elements
of the prime field are existentially \mbox{$\emptyset$-definable}. On the other
hand, many finitely generated extensions of $\Q$ contain existentially
\mbox{$\emptyset$-definable} elements which are transcendental over $\Q$.
Finally, we show that all transcendental elements in $\R$ having a recursive approximation by rationals,
are definable in ${\R}(t)$, and the same holds when one replaces $\R$ by any Pythagorean subfield of $\R$.
\vskip 1.0truecm
\par
\noindent
{\bf 1.~Introduction}
\vskip 0.2truecm
\par
Let ${\cal L}$ be an elementary language. Let ${\cal A}$ be any
\mbox{${\cal L}$-structure} and let $R$ be any \mbox{$n$-ary} relation on $|{\cal A}|$.
Svenonius' theorem (\cite{Svenonius},~[11,~p.~184])
states that the following conditions are equivalent:
\vskip 0.2truecm
\par
$R$ is \mbox{$\emptyset$-definable} in ${\cal A}$ by a formula of ${\cal L}$;
\vskip 0.2truecm
\par
for each elementary extension $({\cal B},S)$ of $({\cal A},R)$
each automorphism $g$ of ${\cal B}$ satisfies
\par
$g(S)=S$.
\vskip 0.2truecm
\par
\noindent
Applying this theorem for fields we conclude that for any field $\K$
and any $r \in \K$ the set $\{r\}$ is \mbox{$\emptyset$-definable} in $\K$
if and only if $g(r)=r$ for each field automorphism $g:\LL \to \LL$
and for each field $\LL$ being an elementary extension of $\K$.
In the next section we give another description of such elements $r$.
\vskip 0.2truecm
\par
\noindent
{\bf 2.~An arithmetic characterization of \mbox{$\emptyset$-definable} elements}
\vskip 0.2truecm
\par
Let $\K$ be a field and let $A$ be a subset of $\K$.
We say that a map $f:A \to \K$ is {\sl arithmetic}
if it satisfies the following conditions:
\par
\noindent
{\bf (1)}~~~if $1 \in A$ then $f(1)=1$,
\par
\noindent
{\bf (2)}~~~if $a,b \in A$ and $a+b \in A$ then $f(a+b)=f(a)+f(b)$,
\par
\noindent
{\bf (3)}~~~if $a,b \in A$ and $a \cdot b \in A$ then $f(a \cdot b)=f(a) \cdot f(b)$.
\par
\noindent
Obviously, if $f:A \to \K$ satisfies condition {\bf (2)}
and $0 \in A$, then $f(0)=0$.
We call an element $r \in \K$ {\sl arithmetically fixed} if there is
a finite set $A(r) \subseteq \K$ (an {\sl arithmetic neighbourhood} of $r$)
with $r \in A(r)$ such that each arithmetic map $f: A(r) \to \K$ fixes~$r$,
i.e. $f(r)=r$. Note that any finite set containing an arithmetic neighbourhood or $r$ is
itself an arithmetic neighbourhood of $r$. We denote the set of arithmetically
fixed elements of a field $\K$ by $\wt{\K}$.
\vskip 0.2truecm
\par
\noindent
{\bf Proposition} (\cite{Tyszka}). $\wt{\K}$ is a subfield of $\K$.
\vskip 0.2truecm
\par
\noindent
{\it Proof.}
We set
$A(0)=\{0\}$ and $A(1)=\{1\}$,
so $0,1 \in \wt{\K}$.
If $r \in \wt{\K}$ then $-r \in \wt{\K}$,
to see this we set
$A(-r)=\{0,-r\} \cup A(r)$.
If $r \in \wt{\K} \setminus \{0\}$
then $r^{-1} \in \wt{\K}$, to see this we set
$A(r^{-1})=\{1,r^{-1}\} \cup A(r)$.
If $r_1,r_2 \in \wt{\K}$ then $r_1+r_2 \in \wt{\K}$,
to see this we set
$A(r_1+r_2)=\{r_1+r_2\} \cup A(r_1) \cup A(r_2)$.
If $r_1,r_2 \in \wt{\K}$ then $r_1 \cdot r_2 \in \wt{\K}$,
to see this we set
$A(r_1 \cdot r_2)=\{r_1 \cdot r_2\} \cup A(r_1) \cup A(r_2)$.
\newline
\rightline{$\Box$}
\par
\noindent
{\bf Theorem 1.} $\wt{\K}=\{x \in \K:~~\{x\}$ is existentially
first-order definable in the language of rings without parameters\}.
\vskip 0.2truecm
\par
\noindent
{\it Proof.} Let $r \in \K$ be arithmetically fixed, and let
$A(r)=\{x_1,...,x_n\}$ be an arithmetic neighbourhood of $r$ with
$x_i \neq x_j$ if $ i \neq j$, and $x_1=r$.
We choose all formulae
$x_i=1$ ($i \in \{1,...,n\}$), $x_i+x_j=x_k$, $x_i \cdot x_j=x_k$ ($i,j,k \in \{1,...,n\}$) that are satisfied
in $A(r)$. Joining these formulae with conjunctions we get
some formula $\Phi$. Let ${\cal V}$ denote the set of variables
in $\Phi$, $x_1 \in {\cal V}$ since otherwise for any $s \in {\K} \setminus \{r\}$
the mapping $f:={\rm id}\left(A(r)\setminus\{r\}\right) \cup \left\{(r,s)\right\}$
satisfies conditions {\bf (1)}-{\bf (3)} and $f(r) \neq r$.
The formula $\underbrace{...~\exists x_i~...}_{\textstyle {x_i \in {\cal V},~i \neq 1}} \Phi$
is satisfied in $\K$ if and only if $x_1=r$. It proves the inclusion $\subseteq$.
We begin the proof of the inclusion $\supseteq$. The
proof presented here is formally a proof by induction on
the complexity of the formula.
We are going to use the following two algorithms.
\newpage
\par
\noindent
{\sc Algorithm 1.} In formulae $\Psi$ of the language of rings,
negations of atomic subformulae are replaced by atomic
formulae. For the language of rings, each negation of an
atomic formula is equivalent to the formula of the form
$W(y_1,...,y_n) \neq~0$, where $y_1$, ..., $y_n$ variables and
$W(y_1,...,y_n) \in {\Z}[y_1,...,y_n]$.
The algorithm selects a variable~$t$ which does not occur in $\Psi$,
and instead of $W(y_1,...,y_n) \neq 0$ introduces to~$\Psi$ the
formula $W(y_1,...,y_n) \cdot t-1=0$. The received formula
features one negation fewer and one variable more.
\vskip 0.2truecm
\par
\noindent
{\sc Algorithm 2.} In formulae $\Psi$ of the language of rings,
some atomic subformulae are replaced by other atomic
formulae or conjunctions of atomic formulae. Atomic
subformulae of the form $y_i+y_j=y_k$, $y_i \cdot y_j=y_k$, $y_i=1$,
($y_i$, $y_j$, $y_k$ variables) are left without changes.
Atomic subformulae of the form $y_i=0$ ($y_i$ is a variable)
are replaced by $y_i+y_i=y_i$. Operation of the algorithm on other
atomic subformulae will be explained on the example of subformula
$1+x+y^2=0$, which is replaced by
$$
(t=1) \wedge (t+x=u) \wedge (y \cdot y=z) \wedge (u+z=s) \wedge (s+s=s),
$$
where variables $t$, $u$, $z$, $s$ do not occur in $\Psi$.
The above conjunction equivalently presents the condition
$1+x+y^2=0$ and is composed solely of the formulae of the form 
$y_i+y_j=y_k$, $y_i \cdot y_j=y_k$, $y_i=1$,
where $y_i$, $y_j$, $y_k$ variables.
\vskip 0.2truecm
\par
\noindent
We start the main part of the proof.
Let $r \in \K$, $\Gamma(x,x_1,...,x_n)$ be a quantifier-free formula
of the language of rings, and
$$
\{r\}=\{x \in \K:~~\K \models \exists x_1 ... \exists x_n~\Gamma(x,x_1,...,x_n)\}
$$
We may assume that $\Gamma(x,x_1,...,x_n)$ has the form
$\Lambda_1 \vee ... \vee \Lambda_l$, where each of the formulae
$\Lambda_1$, ..., $\Lambda_l$ is the conjunction of atomic formulae
and negations of atomic formulae. We want to prove that $r \in \wt{\K}$.
After an iterative application of algorithm~1 to the formula
$\Gamma(x,x_1,...,x_n)$ we receive a quantifier-free formula
$\Omega(x,x_1,...,x_m)$ for which: $\Omega(x,x_1,...,x_m)$ has the form
$\Xi_1 \vee ... \vee \Xi_l$, and each of the formulae
$\Xi_1$, ..., $\Xi_l$ is the conjunction of atomic formulae, and
$$
\{r\}=\{x \in \K:~~ \K \models \exists x_1 ... \exists x_m~\Omega(x,x_1,...,x_m)\},
$$
where $m-n$ is the number of negations in the formula
$\Gamma(x,x_1,...,x_n)$. After an iterative application
of algorithm~2 to the formula $\Omega(x,x_1,...,x_m)$
we receive a quantifier-free formula $\Delta(x,x_1,...,x_p)$
for which: $\Delta(x,x_1,...,x_p)$ has the form
$\Pi_1 \vee ... \vee \Pi_l$, and each of the formulae
$\Pi_1$, ..., $\Pi_l$ is the conjunction of atomic formulae of the form
$y_i+y_j=y_k$, $y_i \cdot y_j= y_k$, $y_i= 1$,
where $y_i$, $y_j$, $y_k$ variables, and
$$
\{r\}=\{x \in \K:~ \K \models \exists x_1 ... \exists x_p ~\Delta(x,x_1,...,x_p)\}
$$
Since
\begin{eqnarray*}
\{r\}=\{x \in \K:~~ \K \models \exists x_1 ... \exists x_p
~\Delta(x,x_1,...,x_p)\}=\\
\bigcup_{i=1}^l \{x \in \K:~
\K \models \underbrace{...~\exists x_s~...}_{x_s \in {\rm Fr}(\Pi_i)\setminus\{x\}}\Pi_i(x,~...,~x_s,~...)\},
\end{eqnarray*}
for some $i \in \{1,...,l\}$ the condition
$$
\{r\}=\{x \in \K:~
\K \models \underbrace{...~\exists x_s~...}_{x_s \in {\rm Fr}(\Pi_i)\setminus\{x\}}\Pi_i(x,~...,~x_s,~...)\}
$$
is satisfied. For indices $s$ for which $x_s$ is a variable in $\Pi_i$,
we choose $w_s \in \K$ for which
$\K \models \Pi_i[x \to r,~...,~x_s \to w_s,~...]$.
Then $A(r)=\{1,r,...,w_s,...\}$ is an arithmetic neighbourhood of
$r$, so $r \in \wt{\K}$.
\newline
\rightline{$\Box$}
\vskip 0.2truecm
\par
Let $\K$ be a field extending $\Q$.
R.~M.~Robinson proved in \cite{Robinson1951}:
if each element of $\K$ is algebraic over $\Q$
and $r \in \K$ is fixed for all automorphisms of $\K$,
then there exist $U(y),V(y) \in {\Q}[y]$ such that $\{r\}$ is
definable in $\K$ by the formula
\vskip 0.2truecm
\par
\noindent
\centerline{$\exists y \left(U(y)=0 \wedge x=V(y)\right)$}
\vskip 0.2truecm
\par
\noindent
{\bf Corollary 1.} If a field $\K$ extends $\Q$ and
each element of $\K$ is algebraic over $\Q$, then
$$
\wt{\K}=\bigcap_{\textstyle \sigma\in{\rm Aut}(\K)}^{}\limits \{x \in \K:~ \sigma(x)=x \}
$$
\par
For a more general theorem and its proof, see [9,~Proposition~1].
Let ${\R}^{\rm alg}:=\{x \in \R:~ x {\rm ~is~algebraic~over~} \Q\}$
and ${\Q}_p^{\rm alg}:=\{x \in {\Q}_p:~ x {\rm ~is~algebraic~over~} \Q\}$.
By Corollary~1, $\wt{{\R}^{\rm alg}}={\R}^{\rm alg}$
and $\wt{{\Q}_p^{\rm alg}}={\Q}_p^{\rm alg}$. It gives
$\wt{\R}={\R}^{\rm alg}$ and $\wt{{\Q}_p}={\Q}_p^{\rm alg}$, see \cite{Tyszka}.
\vskip 0.2truecm
\par
\noindent
{\bf Theorem 2.} Let $\K$ be a field extending $\Q$,
$\phi(x,x_1,...,x_n)$ is a quantifier-free formula of the language of rings, and
$\K \models \exists x \exists x_1 ... \exists x_n \phi(x,x_1,...,x_n)$.
Then there exist a prime number $p$ and $U(y),V(y) \in {\Q}[y]$ such that
\vskip 0.2truecm
\par
\noindent
\centerline{$\{x \in {\Q}_p:~~ {\Q}_p \models \exists x_1 ... \exists x_n \exists y \left( \phi(x,x_1,...,x_n) \wedge U(y)=0 \wedge x=V(y) \right)\}=\{b\}$}
\vskip 0.2truecm
\par
\noindent
for some $b \in {\Q}_p^{\rm alg}$.
\vskip 0.2truecm
\par
\noindent
{\it Proof.} We choose $a,a_1,...,a_n \in \K$ such that
$\K \models \phi(a,a_1,...,a_n)$,
so ${\Q}(a,a_1,...,a_n) \models \exists x \exists x_1 ... \exists x_n \phi(x,x_1,...,x_n)$.
There is a prime number $p$ such that ${\Q}(a,a_1,...,a_n)$ embeds in
${\Q}_p$, see Theorem 1.1 in chapter 5 of \cite{Cassels}.
By this, ${\Q}_p \models \exists x \exists x_1 ... \exists x_n \phi(x,x_1,...,x_n)$.
Since ${\Q}_p^{\rm alg}$ is an elementary subfield of ${\Q}_p$ (\cite{MacIntyre}),
there exists $b \in {\Q}_p^{\rm alg}$ such that
${\Q}_p^{\rm alg} \models \exists x_1 ... \exists x_n \phi(b,x_1,...,x_n)$.
By Robinson's theorem there exist $U(y),V(y) \in {\Q}[y]$ such that $\{b\}$ is
definable in ${\Q}_p^{\rm alg}$ by the formula $\exists y \left(U(y)=0 \wedge x=V(y)\right)$.
Thus,
\begin{eqnarray*}
\{x \in {\Q}_p:~~ {\Q}_p \models \exists x_1 ... \exists x_n \exists y \left( \phi(x,x_1,...,x_n) \wedge U(y)=0 \wedge x=V(y) \right)\}
&=&\\
\{x \in {\Q}_p^{\rm alg}:~~ {\Q}_p^{\rm alg} \models \exists x_1 ... \exists x_n \exists y \left( \phi(x,x_1,...,x_n) \wedge U(y)=0 \wedge x=V(y) \right)\}
&=& \{b\}
\end{eqnarray*}
\rightline{$\Box$}
\vskip 0.2truecm
\par
\noindent
{\bf 3.~Fields with algebraically closed subfields}
\vskip 0.2truecm
\par
We use below "bar" to denote the algebraic closure of a field.
It was proved in~\cite{Tyszka} that $\wt{\C}=\Q$.
Similarly, $\wt{\overline{\Q}}=\Q$.
\vskip 0.2truecm
\par
\noindent
{\bf Theorem 3.} If $\K$ is a field and some subfield of $\K$
is algebraically closed, then $\wt{\K}$ is the prime field in $\K$.
\vskip 0.2truecm
\par
\noindent
{\it Proof.} For any field $\K$ of non-zero characteristic
$\wt{\K}$ is the prime field in $\K$, see \cite{Tyszka}.
Let ${\rm char}(\K)=0$. We may assume that $\K$ extends $\Q$.
By the assumption of the theorem $\K$ extends $\overline{\Q}$.
By the Proposition $\wt{\K} \supseteq \Q$.
We want to prove $\wt{\K} \subseteq \Q$ in a constructive way
without the use of Theorem 1.
Let $r \in \wt{\K}$, and let
$A(r)=\{x_1,...,x_n\}$ be an arithmetic neighborhood of
$r$, $x_i \neq x_j$ if $ i \neq j$,
and $x_1=r$.
We choose all formulae
$x_i=1$ ($i \in \{1,...,n\}$), $x_i+x_j=x_k$, $x_i \cdot x_j=x_k$ ($i,j,k \in \{1,...,n\}$) that are satisfied
in $A(r)$. Joining these formulae with conjunctions we get
some formula~$\Phi$.
Let ${\cal V}$ denote the set of variables
in $\Phi$, $x_1 \in {\cal V}$ since otherwise for any $s \in {\K} \setminus \{r\}$
the mapping $f:={\rm id}\left(A(r)\setminus\{r\}\right) \cup \left\{(r,s)\right\}$
satisfies conditions {\bf (1)}-{\bf (3)} and $f(r) \neq r$. Since $A(r)$ is an arithmetic
neighbourhood of $r$,
\begin{equation}
\tag*{{\bf (4)}}
{\rm ~the~formula~}\underbrace{...~\exists x_i~...}_{\textstyle {x_i \in {\cal V},~i \neq 1}} \Phi
{\rm ~is~satisfied~in~} \K {\rm ~if~and~only~if~} x_1=r
\end{equation}
\vskip 0.2truecm
\par
\noindent
Since $\overline{\K}$ extends $\K$,
$$ \overline{\K} \models
\underbrace{...~\exists x_i~...}_{\textstyle {x_i \in {\cal V},~i \neq 1}} \Phi
[x_1 \to r]
$$
$\overline{\Q}$ is an elementary subfield of $\overline{\K}$ ([6,~p.~306]),
so there exists $r_1 \in \overline{\Q}$ satisfying
\begin{equation}
\tag*{{\bf (5)}}
\overline{\Q} \models
\underbrace{...~\exists x_i~...}_{\textstyle {x_i \in {\cal V},~i \neq 1}} \Phi
[x_1 \to r_1]
\end{equation}
$\K$ extends $\overline{\Q}$, so by {\bf (4)}
there is a unique $r_1 \in \overline{\Q}$ satisfying {\bf (5)} and
this $r_1$ equals $r$. Thus, $r \in \overline{\Q}$ and the formula
$\underbrace{...~\exists x_i~...}_{\textstyle {x_i \in {\cal V},~i \neq 1}} \Phi$
is satisfied in $\overline{\Q}$ if and only if $x_1=r$.
Hence $r \in \wt{\overline{\Q}}=\Q$.
\newline
\rightline{$\Box$}
\vskip 0.2truecm
\par
\noindent
{\bf Corollary 2.} Let $\K$ be an arbitrary field.
Then no subfield of $\wt{\K}$ is algebraically closed.
\vskip 0.2truecm
\par
\noindent
{\bf Theorem 4.} If a field $\K$ extends $\Q$ and $r \in \wt{\K}$, then
$\{r\}$ is definable in $\K$ by a formula
of the form $\exists x_1 ... \exists x_m T(x,x_1,...,x_m)=0$,
where $m \in \{1,2,3,...\}$ and $T(x,x_1,...,x_m) \in {\Z}[x,x_1,...,x_m]$.
\vskip 0.2truecm
\par
\noindent
{\it Proof.} From the definition of $\wt{\K}$ it follows
that $\{r\}$ is definable in $\K$ by a finite system~{\bf (S)} of polynomial
equations of the form $x_i+x_j-x_k=0$, $x_i \cdot x_j-x_k=0$, $x_i-1=0$,
cf. the proof of the inclusion $\subseteq$ inside the proof of Theorem 1.
If~~$\overline{\Q} \subseteq \K$, then by Theorem~3 each element of $\wt{\K}$
is definable in $\K$ by a single equation $w_1 \cdot x+w_0=0$, where
$w_0 \in \Z$, $w_1 \in \Z \setminus \{0\}$.
If~~$\overline{\Q} \not\subseteq \K$, then there exists a polynomial
$$
a_{n}x^{n}+a_{n-1}x^{n-1}+...+a_{1}x+a_{0} \in {\Z}[x]~~~~(n \geq 2,~a_n \neq 0)
$$
having no root in $\K$. By this, the polynomial
$$
B(x,y):=a_{n}x^{n}+a_{n-1}x^{n-1}y+...+a_{1}xy^{n-1}+a_{0}y^{n}
$$
satisfies
\begin{equation}
\tag*{{\bf (6)}}
\forall u,v \in \K \left((u=0 \wedge v=0) \Longleftrightarrow B(u,v)=0 \right),
\end{equation}
see [3,~pp.~363--364] and [14,~p.~108], cf.~[2,~p.~172].
Applying {\bf (6)} to {\bf (S)} we obtain
that {\bf (S)} is equivalent to a single equation
$T(x,x_1,...,x_m)=0$, where $m \in \{1,2,3,...\}$ and $T(x,x_1,...,x_m) \in {\Z}(x,x_1,...,x_m)$.
\newline
\rightline{$\Box$}
\vskip 0.2truecm
\par
Theorem 4 remains true if ${\rm char}(\K)=p \neq 0$.
In this case $\wt{\K}$ is the prime field in~$\K$~(\cite{Tyszka}),
so each element of $\wt{\K}$ is definable by the equation
$w_1 \cdot x+w_0=0$ for some $w_0 \in \{0,1,...,p-1\}$,
$w_1 \in \{1,...,p-1\}$.
\vskip 0.2truecm
\par
\noindent
{\bf 4.~Transcendental elements in finitely generated fields}
\vskip 0.2truecm
\par
It is known (\cite{Koenigsmann}) that for any field $\K$
there is a function field ${\F}/{\K}$ in one variable containing
elements that are transcendental over $\K$ and
first-order definable in the language of rings with parameters from $\K$.
We present similar results with quite different proofs.
\vskip 0.2truecm
\par
\noindent
{\bf Theorem 5.} Let $w$ be transcendental over $\Q$ and
a field $\K$ be finitely generated over ${\Q}(w)$.
Let $g(x,y) \in {\Q}[x,y]$, there exists $z \in \K$ with $g(w,z)=0$,
and the equation $g(x,y)=0$ defines an irreducible algebraic curve of
genus greater than $1$. We claim that some element of $\wt{\K}$ is
transcendental over $\Q$.
\vskip 0.2truecm
\par
\noindent
{\it Proof.} By Faltings' finiteness theorem (\cite{Faltings},~cf.~[8, p.~12], formerly Mordell's conjecture)
the set
$$P:=\{u \in {\K}:~~\exists s \in {\K}~g(u,s)=0\}$$
is finite, $w \in P$. Let $P=\{u_1,...,u_n\}$, $u_i \neq u_j$ if $i \neq j$,
and
$$t_k(x_1,...,x_n):=\sum_{1 \leq i_1<i_2<...<i_k \leq n}\limits
x_{i_1} x_{i_2}...x_{i_k}~~~~~~(k \in \{1,...,n\})$$
denote the basic symmetric polynomials. We claim
that
\begin{equation}
\tag*{{\bf (7)}}
t_1(u_1,...,u_n),...,t_n(u_1,...,u_n) \in \wt{\K}
\end{equation}
and $t_i(u_1,...,u_n)$ is transcendental over $\Q$ for some $i \in \{1,...,n\}$.
We want to prove~{\bf (7)} in a constructive way without the use of Theorem 1.
To prove {\bf (7)} we choose $z_k \in \K$ ($k \in \{1,...,n\}$)
that satisfy $g(u_k,z_k)=0$. There exist $m \in \{1,2,3,...\}$
and
$$
h: \{0,...,m\} \times \{0,...,m\} \to W(m):= \{0\} \cup \left\{\frac{\textstyle c}{\textstyle d}:~~c,d \in \{-m,...,-1,1,...,m\}\right\}
$$
such that
$$g(x,y)=\sum_{i,j \in \{0,...,m\}} h(i,j) \cdot x^{i} \cdot y^{j}$$
Let
\vskip 0.2truecm
\par
\noindent
$$M_{k}:=\left\{u_{i_1}u_{i_2}...u_{i_k}:~~1 \leq i_1<i_2<...<i_k \leq n \right\}~~~~(k \in \{1,...,n\})$$
\vskip 0.2truecm
\par
\noindent
$$N:=\left\{b \cdot u_k^i \cdot z_k^j:~~b \in W(m),~i,j \in \{0,...,m\},~k \in \{1,...,n\} \right\}$$
\vskip 0.2truecm
\par
\noindent
$$
T:=\left\{\sum_{a \in S}\limits a:~~\emptyset \neq S \subseteq N \cup \bigcup_{k=1}^n M_k \right\}
\cup
\left\{u_i-u_j,~\frac{1}{u_i-u_j}:~~i,j \in \{1,...,n\},~i \neq j\right\}
$$
\vskip 0.2truecm
\par
\noindent
Since $M_k \subseteq T$ for each $k \in \{1,...,n\}$,
$t_k(u_1,...,u_n)=\sum_{a \in M_{\scriptstyle k}}\limits a \in T$
for each $k \in \{1,...,n\}$. We claim that $T$ is an arithmetic
neighbourhood of $t_k(u_1,...,u_n)$ for each $k \in \{1,...,n\}$. To prove it
assume that $f:T \to \K$ satisfies conditions {\bf (1)}-{\bf (3)}.
Since $T \supseteq N \supseteq W(m)$, $f$ is the identity on $W(m)$.
For any $k \in \{1,...,n\}$ and any non-empty
$L \varsubsetneq \{0,...,m\} \times \{0,...,m\}$ the elements
$\sum_{(i,j) \in L}\limits h(i,j) \cdot u_k^i \cdot z_k^j$ and
$\sum_{(i,j) \in \left(\{0,...,m\} \times \{0,...,m\} \right) \setminus L}\limits h(i,j) \cdot u_k^i \cdot z_k^j$
belong to $T$. By these facts and by induction
\begin{eqnarray*}
0=f(0)=f(g(u_k,z_k))=
f\left(\sum_{i,j \in \{0,...,m\}} h(i,j) \cdot u_k^i \cdot z_k^j\right)
&=&\\
\sum_{i,j \in \{0,...,m\}} f\left(h(i,j) \cdot u_k^i \cdot z_k^j \right)
=
\sum_{i,j \in \{0,...,m\}} h(i,j) \cdot f(u_k)^{i} \cdot f(z_k)^{j}
=g(f(u_k),f(z_k))
\end{eqnarray*}
for any $k \in \{1,...,n\}$. Thus, $f(u_k) \in P$ for each $k \in \{1,...,n\}$.
Since $1=f(1)=f\left((u_k-u_l) \cdot \frac{1}{u_k-u_l}\right)=\left(f(u_k)-f(u_l)\right) \cdot f\left(\frac{1}{u_k-u_l}\right)$,
we conclude that $f(u_k) \neq f(u_l)$ if $k \neq l$.
Therefore, $f$ permutes the elements of $\{u_1,...,u_n\}$.
By this,
\begin{eqnarray*}
t_k(u_1,...,u_n)=t_k(f(u_1),...,f(u_n))=
\sum_{1 \leq i_1<i_2<...<i_k \leq n} f(u_{i_1})f(u_{i_2})...f(u_{i_k})&=&\\ \\
\sum_{1 \leq i_1<i_2<...<i_k \leq n} f(u_{i_1} u_{i_2} ... u_{i_k})=
f\left(\sum_{1 \leq i_1<i_2<...<i_k \leq n}\limits u_{i_1} u_{i_2} ...u_{i_k}\right)=
f(t_k(u_1,...,u_n))
\end{eqnarray*}
for any $k \in \{1,...,n\}$. We have proved that $T$ is an arithmetic neighbourhood of
$t_k(u_1,...,u_n)$ for each $k \in \{1,...,n\}$, so
$t_k(u_1,...,u_n) \in \wt{\K}$ for each $k \in \{1,...,n\}$.
\vskip 0.2truecm
\par
We prove now that $t_i(u_1,...,u_n)$ is transcendental over $\Q$ for some
$i \in \{1,...,n\}$.
Assume, on the contrary, that all $t_k(u_1,...,u_n)$~($k \in \{1,...,n\}$)
are algebraic over $\Q$. Since $u_1,...,u_n$ are the roots of the polynomial
$$
x^n-t_1(u_1,...,u_n)x^{n-1}+t_2(u_1,...,u_n)x^{n-2}-...+(-1)^{n}t_n(u_1,...,u_n),
$$
we conclude that $u_1,...,u_n$ are also algebraic over $\Q$.
It is impossible, because among elements $u_1,...,u_n$
is $w$ that is transcendental over $\Q$.
\newline
\rightline{$\Box$}
\vskip 0.2truecm
\par
In the proof of Theorem 5 for each $k \in \{1,...,n\}$ the
set $\{t_k(u_1,...,u_n)\}$ is existentially \mbox{$\emptyset$-definable}
in $\K$ by the formula
\begin{flushright}
$\exists u_1 \exists s_1 ... \exists u_n \exists s_n$~~~~~~~~~~
\end{flushright}
\begin{equation}
\tag*{{\bf (8)}}
(g(u_1,s_1)=0 \wedge ... \wedge g(u_n,s_n)=0~ \wedge~ \underbrace{...\wedge u_i \neq u_j \wedge...}_{1 \leq i<j \leq n} ~\wedge ~v=t_k(u_1,...,u_n))
\end{equation}
Applying Theorem 1 we obtain $t_k(u_1,...,u_n) \in \wt{\K}$ for
each $k \in \{1,...,n\}$, unfortunately, without a direct description
of any arithmetic neighbourhood of $t_k(u_1,...,u_n)$. This gives a
non-constructive proof of Theorem 5.
\vskip 0.2truecm
\par
Formula {\bf (8)} has a form
$\exists u_1 \exists s_1 ... \exists u_n \exists s_n \phi(v,u_1,s_1,...,u_n,s_n)$,
where $\phi(v,u_1,s_1,...,u_n,s_n)$ is quantifier-free. By Theorem 2 there exist
a prime number $p$ and $U(y),V(y) \in {\Q}[y]$ such that the formula
\vskip 0.2truecm
\par
\noindent
\centerline{$\exists u_1 \exists s_1 ... \exists u_n \exists s_n \exists y \left( \phi(v,u_1,s_1,...,u_n,s_n) \wedge U(y)=0 \wedge v=V(y) \right)$}
\vskip 0.2truecm
\par
\noindent
defines in ${\Q}_p$ an element that is algebraic over $\Q$.
\vskip 0.2truecm
\par
The proof of Theorem 5 gives an element of $\wt{\K}$
that is transcendental over~$\Q$. Let $\K$ be a field extending $\Q$ and $v \in \wt{\K}$
is transcendental over $\Q$. Since $\wt{\K}$ is a subfield of $\K$, ${\Q}(v) \setminus \Q \subseteq \wt{\K}$.
Obviously, each element of ${\Q}(v) \setminus \Q$ is transcendental over $\Q$.
\vskip 0.2truecm
\par
There exists a function field $\K/\Q$ in one variable such that
$$
\wt{\K}=\K \varsupsetneq \Q=\{x \in \K:~ x {\rm ~is~algebraic~over~} \Q\}
$$
It follows from Proposition 3 in \cite{Lettl}.
\vskip 0.2truecm
\par
Theorem 5 admits a more general form. Let the fields $\K$ and
$\LL$ be finitely generated over $\Q$ such that $\LL$ extends $\K$.
Let $w \in \LL$ be transcendental over $\K$, $g(x,y) \in {\Q}[x,y]$,
there exists $z \in \LL$ with $g(w,z)=0$, and the equation $g(x,y)=0$
defines an irreducible algebraic curve of genus greater than $1$.
Analogously as in the proof of Theorem 5 we conclude that there is
an element of $\wt{\LL}$ that is transcendental over $\K$.
\vskip 0.2truecm
\par
Let $p$ be a prime number, ${\R}(x,y)$~(${\Q}_p(x,y)$)
denote the function field defined by $px^4+p^2y^4=-1$.
The genus of the extension
${\R}(x,y)/{\R}$~(${\Q}_p(x,y)/{\Q}_p$)
is greater than~$1$. By the results in
[7, p.~952, item 3 inside the proof of Theorem 1]
the sets
$$
\{(u,v) \in {\R}(x,y) \times {\R}(x,y):~~pu^4+p^2v^4=-1\}
\setminus
\{(u,v) \in \R \times \R:~~pu^4+p^2v^4=-1\}
$$
$$
\{(u,v) \in {\Q}_p(x,y) \times {\Q}_p(x,y):~~pu^4+p^2v^4=-1\}
\setminus
\{(u,v) \in {\Q}_p \times {\Q}_p:~~pu^4+p^2v^4=-1\}
$$
are finite. Since
$$
\{(u,v) \in \R \times \R:~~pu^4+p^2v^4=-1\}=\emptyset
$$
and
$$
\{(u,v) \in {\Q}_p \times {\Q}_p:~~pu^4+p^2v^4=-1\}=\emptyset,
$$
the sets
$$
\{(u,v) \in {\R}(x,y) \times {\R}(x,y):~~pu^4+p^2v^4=-1\}
$$
and
$$
\{(u,v) \in {\Q}_p(x,y) \times {\Q}_p(x,y):~~pu^4+p^2v^4=-1\}
$$
are finite. Analogously as in the proof of Theorem 5 we conclude
that there is an element of~~$\wt{\R(x,y)}$~($\wt{{\Q}_p(x,y)}$)
that is transcendental over $\R$~(${\Q}_p$).
\vskip 0.2truecm
\par
\noindent
{\bf 5. Recursively defined transcendentals in function fields over archimedean pythagorean fields}
\vskip 0.2truecm
\par
A real number $r$ is called recursively approximable,
if there exists a computable sequence of
rational numbers which converges to $r$, see \cite{Zheng}.
Let $\omega:=\{0,1,2,...\}$, $\K$ be a subfield of $\R$.
$\K$ is said to be Pythagorean if
$$
\forall x \in \K ~(0 \leq x \Rightarrow \exists y \in \K ~~x=y^2)
$$
Our next theorem is inspired by Cherlin's example in [7,~p.~949].
\vskip 0.2truecm
\par
\noindent
{\bf Theorem 6.} If $\K$ is a Pythagorean subfield of $\R$,
$t$ is transcendental over $\K$, and
$r \in \K$ is recursively approximable,
then $\{r\}$ is \mbox{$\emptyset$-definable} in $({\K}(t),+,\cdot,0,1)$.
\vskip 0.2truecm
\par
\noindent
{\it Proof.} It follows from [13,~p.~280] that
there is a formula ${\cal N}(x)$ in the language of rings such that
\begin{equation}
\tag*{{\bf (9)}}
\{x \in {\K}(t):~{\K}(t) \models {\cal N}(x)\}=\omega
\end{equation}
Let ${\cal M}(x)$ abbreviate $\exists y ~1+x^4=y^2$.
It is known that
$$
\{x \in {\K(t)}:~{\K}(t) \models {\cal M}(x)\}=\K,
$$
for the proof see [6,~p.~34]. Assume that $r \geq 0$,
the proof in case $r \leq 0$ goes analogically.
There exist recursive functions $f:\omega \to \omega$
and $g:\omega \to \omega \setminus \{0\}$ such that
$\lim_{n \to \infty}\limits \frac{f(n)}{g(n)}=r$.
Since $f$ and $g$ are recursive
there exist formulae $F(s,t)$ and $G(s,t)$ (both in the language of rings)
for which
$$
\forall n,m \in \omega ~(m=f(n) \Longleftrightarrow \omega \models F(n,m))
$$
and
$$
\forall n,m \in \omega ~(m=g(n) \Longleftrightarrow \omega \models G(n,m))
$$
By {\bf (9)} we can find formulae $\wt{F}(s,t)$ and $\wt{G}(s,t)$ for which
$$
\forall s,t \in {\K}(t) ~((s \in \omega \wedge t \in \omega \wedge t=f(s)) \Longleftrightarrow {\K}(t) \models \wt{F}(s,t))
$$
and
$$
\forall s,t \in {\K}(t) ~((s \in \omega \wedge t \in \omega \wedge t=g(s)) \Longleftrightarrow {\K}(t) \models \wt{G}(s,t))
$$
Let $a<b$ abbreviate 
$$
a \neq b \wedge {\cal M}(a) \wedge {\cal M}(b) \wedge \exists c ~({\cal M}(c) \wedge a+c^2=b)
$$
The formula
\begin{eqnarray*}
{\cal M}(x) \wedge \forall \varepsilon ~(0<\varepsilon \Rightarrow \exists z \exists s \exists u \exists v
~(z \neq x ~\wedge~ x<z+\varepsilon ~\wedge~ z<x+\varepsilon ~\wedge\\
{\cal N}(s) \wedge {\cal N}(u) \wedge {\cal N}(v) \wedge
\wt{F}(s,u) \wedge \wt{G}(s,v) \wedge z \cdot v=u))
\end{eqnarray*}
defines $r$ in $({\K}(t),+,\cdot,0,1)$.
\newline
\rightline{$\Box$}
\vskip 0.2truecm
\par
Let ${\cal L}$ be an elementary language, let $M$ be an \mbox{${\cal L}$-structure},
and let $U$ be an \mbox{$n$-ary} relation on~$M$.
We say that $U$ is implicitly \mbox{$\emptyset$-definable}
in $M$ if there exists a sentence~$\Phi$ in the language
${\cal L} \cup \{{\cal U}\}$ with an additional \mbox{$n$-ary}
predicate symbol ${\cal U}$, such that for all \mbox{$n$-ary}
relations $U^{\ast}$ on $M$, $(M,U^{\ast}) \models \Phi$
if and only if $U^{\ast}=U$, see the introductory part of~\cite{Fukuzaki}.
\vskip 0.2truecm
\par
\noindent
{\bf Theorem 7.} If a real number $r$ is recursively approximable,
then $\{r\}$ is existentially \mbox{$\emptyset$-definable}
in $(\R,+,\cdot,0,1,U)$ for some unary predicate $U$ which
is implicitly \mbox{$\emptyset$-definable} in $(\R,+,\cdot,0,1)$.
\vskip 0.2truecm
\par
\noindent
{\it Proof.} If $r$ is a rational number then $\{r\}$
is existentially \mbox{$\emptyset$-definable} in $(\R,+,\cdot,0,1)$.
At this moment we assume that $r$ is an irrational number.
We may assume without loss of generality that $r<0$,
so there exists an integer $i<r$.
There exist recursive functions
$f:\omega \to \omega$ and $g:\omega \to \omega \setminus \{0\}$
such that $\lim_{n \to \infty}\limits -~\frac{f(n)}{g(n)}=r$,
we may assume without loss of generality that
$-~\frac{f(n)}{g(n)} \in (i,0)$ for each $n \in \omega$.
Since $f$ and $g$ are recursive,
there exist formulae $F(s,t)$ and $G(s,t)$
(both in the language of rings) for which
$$
\forall n,m \in \omega
~(m=f(n) \Longleftrightarrow \omega \models F(n,m))
$$
and
$$
\forall n,m \in \omega
~(m=g(n) \Longleftrightarrow \omega \models G(n,m))
$$
Let
$$
U:=\{r+i\} \cup \bigl\{ -~\frac{f(n)}{g(n)}:~~n \in \omega \bigr\} \cup \omega
$$
and ${\cal U}$ be a unary predicate symbol for membership in $U$.
Let $x \leq y$ abbreviate
$$
\exists s ~x+s^2=y,
$$
$x<y$ abbreviate
$$
x \leq y \wedge x \neq y,
$$
${\rm succ}(x,y)$ abbreviate
$$
x<y \wedge {\cal U}(x) \wedge {\cal U}(y) \wedge
\forall z ~((x<z \wedge z<y) \Rightarrow \neg {\cal U}(z)),
$$
${\rm accum}(x)$ abbreviate
$$
\forall \varepsilon ~(0<\varepsilon \Rightarrow \exists z 
~(z \neq x \wedge x<z+\varepsilon \wedge z<x+\varepsilon \wedge {\cal U}(z)))
$$
We have:
$$
\forall x \in \R ~(x \in \omega \Longleftrightarrow \R \models (0 \leq x \wedge {\cal U}(x)))
$$
Therefore, extending the language of rings with predicate symbol
${\cal U}$ for membership in $U$ we can find
formulae $\wt{F}(s,t)$ and $\wt{G}(s,t)$ for which
$$
\forall s,t \in \R
~((s \in \omega \wedge t \in \omega \wedge t=f(s))
\Longleftrightarrow \R \models \wt{F}(s,t))
$$
and
$$
\forall s,t \in \R
~((s \in \omega \wedge t \in \omega \wedge t=g(s))
\Longleftrightarrow \R \models \wt{G}(s,t))
$$
The sentence
\begin{eqnarray*}
{\cal U}(0) \wedge \forall x ~((0 \leq x \wedge {\cal U}(x)) \Rightarrow {\rm succ}(x,x+1))
&\wedge&\\
\forall x ~((0 \leq x+\underbrace{1+...+1}_{|i|{\rm -times}} ~\wedge~ x<0) \Longleftrightarrow
\exists s \exists u \exists v \\
~(0 \leq s \wedge {\cal U}(s) \wedge 0 \leq u \wedge {\cal U}(u) \wedge 0 \leq v \wedge {\cal U}(v)
\wedge \wt{F}(s,u) \wedge \wt{G}(s,v) \wedge u+x \cdot v=0))
&\wedge&\\
\forall x ~((0<x+\underbrace{1+...+1}_{2|i|{\rm -times}} ~\wedge~ x+\underbrace{1+...+1}_{|i|{\rm -times}} < 0)
\Rightarrow ({\cal U}(x) \Longleftrightarrow {\rm accum}(x+\underbrace{1+...+1}_{|i|{\rm -times}})))
&\wedge&\\
\forall x ~(x+\underbrace{1+...+1}_{2|i|{\rm -times}} \leq 0 \Rightarrow \neg {\cal U}(x))
\end{eqnarray*}
is valid in $\R$ if and only if ${\cal U}(x)$ means $x \in U$,
so $U$ is implicitly \mbox{$\emptyset$-definable} in $\R$.
The formula
$$
\exists t \exists y ~(x+t^2=0 \wedge x=y+\underbrace{1+...+1}_{|i|{\rm -times}} ~\wedge~ {\cal U}(y))
$$
defines $r$ in $(\R,+,\cdot,0,1,U)$.
\newline
\rightline{$\Box$}
\vskip 0.2truecm
\par
\noindent
{\bf Acknowledgement.} The author wishes to thank the anonymous referee for valuable
suggestions. The author also thanks the anonymous referee of Bull. Lond. Math. Soc. for
valuable suggestions to the previous, shorter version of the paper considered in 2005.

\newpage
\par
\noindent
Apoloniusz Tyszka\\
Technical Faculty\\
Hugo Ko\l{}\l{}\k{a}taj University\\
Balicka 104, 30-149 Krak\'ow, Poland\\
E-mail address: {\it rttyszka@cyf-kr.edu.pl}\\
\end{document}